\documentclass[11pt]{article}
\usepackage{a4wide}
\usepackage{amsmath}

\usepackage{amsfonts}
\usepackage{amssymb}
\usepackage{euscript}
\newenvironment{proof}{\noindent {\it Proof.~~}\ }{\  \rule{1mm}{2mm}\medskip}

\newenvironment{proof*}{\noindent {\it Proof.~~}\ }{}

\newtheorem{theorem}{Theorem}
\newtheorem{lemma}[theorem]{Lemma}
\newtheorem{corollary}[theorem]{Corollary}
\newtheorem{proposition}[theorem]{Proposition}
\newtheorem{conj}[theorem]{Conjecture}

\def\N{\mathbb N}

\begin{document}
\date{}
\title{On  iterated image size for point-symmetric relations }
\author{Yahya Ould Hamidoune\thanks{Universit\'e Pierre et Marie Curie, Paris.
{\tt yha@ccr.jussieu.fr}}}
\maketitle

\begin{abstract}

Let $\Gamma =(V,E)$  be a point-symmetric reflexive  relation and let $v\in V$ such that
 $|\Gamma (v)|$ is finite (and hence $|\Gamma (x)|$ is finite for all $x$
 , by the transitive action of the group of automorphisms).
Let $j\in \N$ be an integer such that $\Gamma ^j(v)\cap \Gamma ^{-}(v)=\{v\}$. Our main result states that
 $$
|\Gamma ^{j}  (v)|\ge | \Gamma ^{j-1}  (v)| + |\Gamma (v)|-1.$$

As an application we have $
|\Gamma ^{j}  (v)| \ge 1+( |\Gamma (v)|-1)j.$ The last result confirms a recent conjecture of Seymour
in the case of vertex-symmetric graphs. Also it gives a short proof for the validity of the Caccetta-H\"aggkvist conjecture for vertex-symmetric graphs and generalizes an additive result of Shepherdson.

 \end{abstract}

\section{Introduction}
  Let $G$ be an abelian group  and let $A,S$ be finite subsets of $G$ with $0\notin S$.
  Shepherdson's generalization of the Cauchy-Davenport Theorem states that $|A\cup (A+S)|\ge |A|+|S|$ if
  $A\cup (A+S)$ contains no subgroup generated by some element of $S$.

 As an application Shepherdson \cite{sheph} proved  that
there are $s_1, \cdots , s_k\in S$ such that $k\leq \lceil \frac{|G|}{|S|}\rceil $ and  $\sum _{1\leq i \leq k} s_i=0,$ if $G$ is finite. The paper of Shepherdson includes thanks to Heilbronn
for suggesting this application together with a mention that Chowla obtained some related zero-sum  results.

Let $D=(V,E)$  be a loopless finite digraph  with minimal outdgree at least $1$. It is well known that $D$ contains a directed cycle. The smallest cardinality  of such a cycle is called the girth of $D$ and will be denoted by $g(D)$.
In 1970 Behzad, Chartrand and  Wall \cite{Behzad}  conjectured that $|V|\ge r(g(D)-1)+1$, if $d^+(x)=d^-(x)=r$ for all $x\in V$.
In 1978,  Caccetta and H\"aggkvist \cite{CH} made the stronger conjecture :
 $$|V|\ge \min (d^+x: x\in V)(g(D)-1)+1.$$

 These conjectures are still largely open, even for the special case $g(D)=4$. The reader may find references and  results about this question in \cite{Bondy}.

 These conjectures were proved by the author for vertex-symmetric digraphs \cite{HEJC}. This result applied to Cayley graphs shows the validity of Shepherdson's zero-sum result for all finite groups. Unfortunately we were not aware at that moment of Shepherdson's result. Our proof  \cite{HEJC} is based on the properties of atoms of a finite digraph and Menger's Theorem.  A description of Cayley graphs on finite Abelian groups
 such that $|V|= r(g(D)-1)+1,$ where $r$ is the outdegree was obtained by the authors of \cite{AOY} using   Kemperman critical pair Theory \cite{Kem}.
 A new proof of the
 Caccetta and H\"aggkvist conjecture for vertex-symmetric digraphs  based
 on  an additive  result of Kemperman \cite{Kemp} and the representation of vertex symmetric digraphs as coset graphs is  given in \cite{Nat}.

More recently  Seymour proposed the following conjecture \cite{Seymour}:

Let $D$ be a loopless digraph and let $r\ge 1$ be an integer. Then there is a vertex $a$ such that  $$|\Gamma (a)\cup \Gamma ^2 (a)\cup \cdots \cup \Gamma ^{g-2}(a)| \ge r(g-2),$$ where $g=g(D).$

The case   $g=4$ of this conjecture is mentioned in \cite{Bondy}.
 Seymour's Conjecture implies the conjecture of Behzad, Chartrand and  Wall.   Seymour's Conjecture  also implies that
 $D$  contains a directed cycle $C$ with $|C|\leq \lceil \frac{|V|-1}{r} \rceil+1$. Notice that   the Caccetta-H\"aggkvist Conjecture states that  $D$  contains a directed cycle $C$ with  $|C|\leq \lceil \frac{|V|}{r} \rceil$.

We shall allow infinite relations. The classical strong connectivity of digraphs needs to be modified
in this case in order to have a good lower bound of the size of the image of a set. Also the presence of loops will simplify the presentation of the connectivity method. Since this convention is unusual in this part of Graph Theory, we shall work with relations. Our terminology will be developed in the next section.

Seymour's conjecture may be formulated as follows :

\begin{conj} \cite{Seymour} Let $\Gamma =(V,E)$ be a finite reflexive  relation and let $j$ be an integer.
Then there is an $x\in V$ such that one of the following conditions holds.
\begin{itemize}
 \item   $|\Gamma ^{j}(x)|\ge 1+j(|\Gamma (x)|-1)$.
 \item   $\Gamma ^{-1}(x)\cap  \Gamma ^{j}(x)\neq \{x\}$.
 \end{itemize}
\end{conj}

Our main result is the following one:

Let $\Gamma =(V,E)$  be a point-symmetric reflexive  relation and let $v\in V$
 such that $|\Gamma (v)|$ is finite. Let $j\ge 1$ be such that $\Gamma ^j(v)\cap \Gamma ^{-}(v)=\{v\}$. Then
 $$|\Gamma ^{j}  (v)|\ge | \Gamma ^{j-1}  (v)| + |\Gamma (v)|-1.$$

This result implies the validity of the above conjectures for vertex-symmetric graphs.

\section{Terminology }
Let $V$ be a set. The diagonal of $V$ is by definition $\Delta _V=\{(x,x) : x\in V\}$.
Let $E \subset V\times V$.  The ordered pair $\Gamma = (V,E)$ will be called  a  {\em relation }.
The relation $\Gamma$ is said to be {\em
reflexive} if $\Delta _V \subset E.$

Let $ a\in V$ and let $A\subset V$. The image of $a$ is by definition  $$\Gamma (a)=\{x:  (a,x)\in E\}.$$ The image of $A$ is by definition $$\Gamma (A)=\bigcup \limits_{x\in A}
\Gamma (x).$$ The cardinality of the image of $x$ will be called the {\it degree} of $x$ and will be denoted by $d(x)$.
 The relation $\Gamma$ will be called {\em regular} with degree $r$ if
the elements of $V$ have the same  degree $r$.
We shall say that  $\Gamma$ is {\em locally finite} if  $d (x)$ is finite for all $x$.
The {\em reverse   } relation of $\Gamma $ is by definition $\Gamma ^-=(V,E^-)$, where
$E^-=\{(x,y) \Big| \   (y,x) \in  E\}.$
The {\em restriction} of $\Gamma =(V,E)$ to a subset $W\subset V$ is defined as the relation
$\Gamma [W]=(W,E\cap (W\times W))$.

Let $\Phi=(W,F)$ be a relation. A function $f : V  \longrightarrow W$ will be called a {\em homomorphism }if for all $x,y\in V$ such that $(x,y)\in E$, we have
$(f(x),f(y))\in F$.

The relation $\Gamma$ will be called {\em point-symmetric}  if for all $x,y\in V$, there is an automorphism $f$
such that $y=f(x)$.
Clearly a point-symmetric relation is regular.

 We identify graphs and their relations. A loopless finite relation will be called a  {\em digraph}.
The reader may replace everywhere  the term "relation" by "graph".  In this case we mention some differences between our terminology (which follows closely the standard notations of Set Theory) and the notations used in some text books of Graph Theory.  We point out that
our graphs are usually
called directed graphs without multiple arcs or digraphs.
Notice that the notion  $\Gamma (a)$ used here and in Set Theory is written $\Gamma ^+(a)$ in
some text books in Graph Theory. Also our notion of degree is called {\it outdegree}. We made the choice of
Set Theory terminology since some parts of this paper could have some interest in Group Theory and Number Theory.

We shall use the composition  $\Gamma _1 \circ\cdots \circ \Gamma _k$ of relations $\Gamma _1,
\cdots , \Gamma _k$ on $V$. If all these relations are equal to $\Gamma$, we shall write

$$\Gamma _1 \circ \cdots \circ \Gamma _k=\Gamma ^k.$$

We shall write $\Gamma ^{0}$ for the identity relation $I_V=(V,\Delta _V)$. Also we shall write $\Gamma ^{-j}$ instead of
$({\Gamma ^{-}})^j.$

\section{Connectivity}

Let $\Gamma =(V,E)$ be a  relation. For $X\subset V$, we shall write
 $$\partial _{\Gamma}(X)= \Gamma (X)\setminus X .$$
When the context is clear the reference to $\Gamma$ will be omitted.

 Let $\Gamma =(V,E)$ be a locally finite  reflexive relation.
The {\em connectivity} of $\Gamma$ is by definition $\kappa  (\Gamma)=|V|-1$, if $E= V\times V.$
Otherwise
\begin{equation}  \label{eq:kappa}
\kappa  (\Gamma)=\min  \{|\partial (X)|\  : \ \ 1\le |X|<\infty \ \mbox{\rm and }\ \Gamma (X)\neq V\}.
\end{equation}

 A subset $X$ achieving the  minimum in  (\ref{eq:kappa}) is called a
{\em fragment} of $\Gamma$. A fragment with minimum cardinality
is called an {\em atom}. The cardinality of an atom of $\Gamma$ will be denoted by $a(\Gamma)$. It is not true that distinct  atoms are always disjoint. But the author proved in \cite{HATOM} that, if $V$ is finite, then
 distinct  atoms of $\Gamma$ are disjoint, or distinct  atoms of $\Gamma ^-$ are disjoint. In \cite{Hdm},
 it was observed that the same methods imply that distinct  atoms of $\Gamma$ are disjoint if $V$ is infinite. One may find
 in \cite{Hast} unified proofs and some applications to Group Theory and Additive Number Theory.

 As a consequence of this result we could obtain :

\begin{proposition} \cite{HATOM, HJCT, Hdm, Hast} {
Let $\Gamma =(V,E)$ be a
 locally-finite point-symmetric  relation with $E\neq V\times V$. Suppose that $V$ is infinite or  that $ a(\Gamma)  \leq
  a(\Gamma ^-) $.
Let $A$ be an atom of $\Gamma$. Then the subrelation $\Gamma [A]$ induced on $A$ is a
point-symmetric relation. Moreover  $|A|\le \kappa (\Gamma)$.
\label{basic}} \end{proposition}
\section{Iterated image size }

\begin{lemma}\label{circl} Let $\Gamma =(V,E)$ be a point-symmetric relation. Then for all $i$, $\Gamma ^i$ is point-symmetric.

\end{lemma}

\begin{proof}
Clearly any automorphism of $\Gamma $ is an automorphism of $\Gamma ^i$.
\end{proof}

\begin{theorem}\label{main}
Let $\Gamma =(V,E)$  be a point-symmetric reflexive locally finite relation and  let $v\in V.$
Let $j\ge 1$ be an integer such that  $\Gamma ^j(v)\cap \Gamma ^{-}(v)=\{v\}$. Then
 $$
|\Gamma ^{j}  (v)|\ge | \Gamma ^{j-1}  (v)| +
|\Gamma (v)|-1.$$
\end{theorem}

\begin{proof}

Set $V_0=\bigcup _{0\le i} \Gamma ^i (v).$
Clearly $\Gamma ^j(v)\subset V_0.$
So we may assume that $\Gamma=\Gamma [V_0]$ and $V=V_0.$

In the finite case this means that we restrict ourselves to the connected component containing $v$.

We shall assume  $j>1$, since  the result is obvious for $j=1.$

With this hypothesis,  clearly we have  $\kappa (\Gamma )\ge 1.$

Clearly  $$E\neq V\times V.$$
 Set $\kappa =\kappa (\Gamma ).$
Let $A$ be an atom of $\Gamma$ containing   $v$.
The proof is by induction on   $|\Gamma (v)|.$

Put $r=|\Gamma (v)|.$
Assume first $$\kappa =r-1.$$ Observe that $\Gamma ^{j}  (v)\ne V$. Then by the definition of $\kappa$, we have $$|\Gamma ^j(v)\setminus \Gamma ^{j-1}  (v)|=|\partial (\Gamma ^{j-1}(v))|  \geq \kappa =r-1.$$ The result holds in this case. So we may assume
\begin{equation}\label{smallcon}
\kappa \leq r-2,
 \end{equation}
and hence $r\ge 3.$  Then
$|A|\ge 2$, since otherwise
$\kappa =|\partial (A)|=r-1$.

{\bf Case 1}. $V$ is infinite or $ a(\Gamma)  \leq
  a(\Gamma ^-) $.

   By Proposition \ref{basic}, $\Gamma [A]$ is point-symmetric (and hence regular) and
   \begin{equation}\label{eqAsmall}
   |A|\le \kappa.
   \end{equation}
  Put $r_0=|\Gamma (v)\cap A|.$

  Put $X= \Gamma ^{j-1}(v)$ and $Y=A\cup X.$
  By the induction hypothesis, we have

  \begin{equation}\label{eqint}
|\partial (X)\cap A|=|\Gamma ^{j}(v)\cap A|-|\Gamma ^{j-1}(v) \cap A|\geq r_0-1.
  \end{equation}

    Let us prove that
  \begin{equation} \label{eqaug}
  \partial (Y)
  \subset (\partial (X)\setminus A)\cup (\partial (A)\setminus \Gamma (v)).
  \end{equation}

  Since $j\ge 2$, we have $\Gamma (v)\subset X$ and hence $\Gamma (v)\cap \partial (X)=\emptyset.$ Then
  (\ref{eqaug}) clearly holds.

  It follows that

   $$\partial (Y)\setminus \partial (X)\subset  \partial (A)\setminus \Gamma (v),$$

    and hence we have
    \begin{eqnarray*}
  |\partial (Y)\setminus \partial (X)| &\le&|\partial (A)|-  |\partial (A)\cap  \Gamma (v)|\\
  &=& \kappa -|\Gamma (v)|+|A\cap \Gamma (v)|\\
  &=&\kappa +r_0-r.
  \end{eqnarray*}

  Hence
  \begin{equation}\label{eqext}
|\partial (Y)\setminus \partial (X)|\le \kappa +r_0-r.
\end{equation}

  Let us show that  $\Gamma  (Y)\neq V.$ This holds obviously if $V$ is infinite.
  So we may assume $V$ finite. In this case we have $|\Gamma ^-(v)|= |\Gamma (v)|.$

   Clearly we have
  $$\Gamma (Y)=\Gamma (X)\cup (A\setminus \Gamma (v))\cup (\partial (Y)\setminus \partial (X)).$$

  It follows using ( \ref{eqAsmall}) and (\ref{eqext}) that

  \begin{eqnarray*}
  |\Gamma (Y)|&\le&|\Gamma (X)|+| A\setminus \Gamma (v)|+| \partial (Y)\setminus \partial (X)|\\
  &\leq &|V\setminus (\Gamma ^-(v)\setminus \{v\})|+|A|-r_0+\kappa+r_0-r\\
  &= &|V|+|A|+\kappa-2r+1\\
  &\le& |V|+2\kappa-2r+1\le |V|-3.
  \end{eqnarray*}

  By the definition of $\kappa$, we have $|\partial (Y)|\geq \kappa.$

  By (\ref{eqint}) and (\ref{eqext}),
  \begin{eqnarray*}
  |\partial (X)|&=& |\partial (X)\cap A|+|\partial (Y)\cap \partial (X)|\\
   &\ge& r_0-1+|\partial (Y)|-|\partial (Y) \setminus \partial (X)|\\
  &\geq &r_0-1+ \kappa -(\kappa +r_0-r)=r-1,
  \end{eqnarray*}
  and the result is proved since $$\partial (X)= \Gamma ^{j}(v) \setminus \Gamma ^{j-1}(v).$$

{\bf Case 2}. $V$ is finite and $ a(\Gamma) >
  a(\Gamma ^-)$.

  The argument used in Case 1, shows that  $
|\Gamma ^{-j}  (v)\setminus \Gamma ^{-(j-1)}  (v)| \ge r-1.$

  By Lemma \ref{circl}, $\Gamma ^j$ is point-symmetric. Since $V$ is finite, $\Gamma ^j$
  and its reverse $\Gamma ^{-j}$ have the same degree.
  Therefore observing that these relations are reflexive
  \begin{eqnarray*}
   r-1&\leq &|\Gamma ^{-j}  (v)|-| \Gamma ^{-(j-1)}  (v)| \\
&=&|\Gamma ^{j}  (v)|-| \Gamma ^{(j-1)}  (v)|.\\
  \end{eqnarray*}

\end{proof}

The next result shows the validity of the conjecture of Seymour mentioned in the introduction in the case of relations with a symmetric group of automorphisms.

\begin{corollary}\label{Seymour} Let $\Gamma =(V,E)$  be a point-symmetric reflexive   relation with
 degree $r$ and let $v\in V$. Let $j\ge 1$ be an integer such that $\Gamma ^j(v)\cap \Gamma ^{-}(v)=\{v\}$
  then $$
|\Gamma ^{j}  (v)| \ge 1+(r-1)j.$$
\end{corollary}

\begin{proof} The proof follows by induction using Theorem \ref{main}

\end{proof}

\begin{corollary}\cite{HEJC}\label{CHVT} Let $\Gamma =(V,E)$  be a point-symmetric digraph with
 degree $r\ge 1$ and put  $g=g(\Gamma) $.
  Then $
|V| \ge 1+r(g-1).$
\end{corollary}

\begin{proof}
Set $ \Phi =(V,E \cup \Delta _V)$.
Let $v\in V$. Clearly we have  $\Phi ^{g-2}(v)\cap \Phi ^{-}(v)=\{v\}$.
 By Corollary  \ref{Seymour}, $|V|-r=|V\setminus (\Phi ^-(a)\setminus \{a\})|\geq |\Phi ^{g-2}(v)|\ge 1+(g-2)r.$
\end{proof}

This result, proved in \cite{HEJC}, shows the validity of the
Caccetta-H\"aggkvist Conjecture for point-symmetric graphs. But the
proof obtained here is much easier.

\begin{corollary}\cite{HEJC} Let $G$ be a group of order $n$
 and let $S\subset G\setminus \{1\}$ with cardinality $=s$. There are elements $s_1, s_2, \cdots ,
 s_k\in S$  such that $k\leq \lceil \frac{n}{s}\rceil$ and
 $\prod \limits_{1\leq i \leq k} s_i=1$.
\end{corollary}

The proof follows by applying Corollary \ref{CHVT} to the Cayley graph defined by $S$
on $G$. In particular the theorem of Shepherdson mentioned in the introduction holds for all finite groups.


\begin{thebibliography}{99}




\bibitem{Behzad}
M. Behzad, G. Chartrand and C.E. Wall, On minimal regular digraphs
with given girth, Fund. Math. 69 (1970), 227-231.
\bibitem{Bondy}
J. A. Bondy, Counting subgraphs: a new approach to the
Caccetta-H\"aggkvist conjecture. Graphs and combinatorics
(Marseille, 1995). Discrete Math. 165/166 (1997), 71-80.
\bibitem{CH}
L. Caccetta and R. H\"aggkvist, On minimal  digraphs with given girth,  Proceedings of the Ninth Southeastern Conference on Combinatorics, Graph Theory, and Computing (Florida Atlantic Univ., Boca Raton, Fla., 1978), (Winnipeg, Man.),  Congress. Numer., XXI, Utilitas Math.  (1978), 181-187.


\bibitem{HATOM} Y.O. Hamidoune, Sur les atomes d'un graphe orient\'e,
 C.R. Acad. Sc. Paris A 284 (1977),   1253-1256.


\bibitem{HJCT} Y.O. Hamidoune,   Quelques probl\`emes de connexit\'e dans les
graphes orient\'es, J. Comb. Theory B 30 (1981), 1-10.


\bibitem{HEJC} Y.O. Hamidoune,   An application  of  connectivity theory in graphs to
      factorizations of elements in groups, Europ. J of Combinatorics 2 (1981), 349-355.

\bibitem{Hdm} Y.O. Hamidoune,   Sur les atomes d'un graphe de Cayley infini,
 Discrete Math., 73 (1989), 297-300.

 \bibitem{Hast} Y.O. Hamidoune,  On small subset product in a group.
Structure Theory of set-addition, Ast\'erisque. no. 258(1999),
xiv-xv, 281-308.

\bibitem{AOY} Y.O. Hamidoune,  A. Llad\'o and O. Serra,
  Vosperian and superconnected Abelian Cayley digraphs,
 Graphs and Combinatorics 7(1991), 143-152.

\bibitem{Nat} Melvyn B. Nathanson, The Caccetta-H\"aggkvist conjecture
and Additive Number Theory, http://arxiv.org/archive/math: eprint arXiv:math/0603469.
\bibitem{Kemp} J.H.B. Kempermann, On complexes in a semigroup, Nederl. Akad. Wetensch. Proc. Ser. A. 59= Indag. Math. 18(1956), 247-254.
\bibitem{Kem} J. H. B. Kemperman, On small sumsets in a  Abelian group,
{ Acta Math.} 103 (1960), 63-88.
\bibitem{Seymour} P. Seymour, Oral communication.
\bibitem{sheph} J. C. Shepherdson, On the addition of elements of a sequence, J. London Math Soc. 22(1947), 85-88.

\end{thebibliography}
\end{document}